\documentclass{rae}
\usepackage{amsmath}\usepackage{amssymb}

\author{Krzysztof Ciesielski%
\thanks{The first author was partially supported by
NSF Cooperative
Research Grant INT-9600548, with its Polish part being
financed by Polish Academy of Science PAN, and 1996/97 West Virginia
University Senate Research Grant.
\endgraf
\hspace{4pt}
Papers authored or co-authored by a Contributing Editor are managed
by a Managing Editor or one of the other Contributing Editors.
\endgraf
\hspace{4pt}
The authors like to thank the referee for many valuable comments and
suggestions.},
Department of Mathematics, West Virginia University,
Morgantown, WV 26506-6310, USA
(KCies@wvnvms.wvnet.edu)\\
Saharon Shelah%
\thanks{This work was supported in part by a grant from
``Basic Research Foundation''
of the Israel Academy of Sciences and Humanities.
Publication 680.
},
Institute of Mathematics,
the Hebrew University of Jerusalem,
91904 Jerusalem, Israel,
and
Department of Mathematics,
Rudgers University,
New Brunswick, NJ 08854, USA}

\title{Uniformly antisymmetric function with bounded range}
\date{}

\MathReviews{26A15, 26A21}

\newcommand{\real}{{\mathbb R}}

\newcommand{\rational}{{\mathbb Q}}

\newcommand{\N}{{\mathbb N}}

\newcommand{\la}{{\langle}}
\newcommand{\ra}{{\rangle}}

\newcommand{\e}{{\emptyset}}
\newcommand{\ep}{{\varepsilon}}

\renewcommand{\H}{{\cal H}}
\renewcommand{\P}{{\cal P}}

\def\proof{\noindent {\sc Proof. }}

\newtheorem{theorem}{Theorem}
\newtheorem{corollary}[theorem]{Corollary}
\newtheorem{proposition}[theorem]{Proposition}
\newtheorem{lemma}[theorem]{Lemma}
\newtheorem{problem}[theorem]{Problem}
\newtheorem{example}[theorem]{Example}
\newtheorem{definition}[theorem]{Definition}
\newtheorem{remark}[theorem]{Remark}
\newtheorem{Fact}[theorem]{Fact}

\newcommand{\thm}[2]{\begin{theorem}\label{#1}{\sl #2}\end{theorem}}

\begin{document}

\maketitle

\begin{abstract}
The goal of this note is to construct
a uniformly antisymmetric function
$f\colon\real\to\real$ with a bounded countable range.
This answers Problem~1(b) of
Ciesielski and Larson~\cite{CL:Unif}. (See also the list of problems in
Thomson~\cite{thomson} and Problem~2(b)
from Ciesielski's survey \cite{56:surv}.)
A problem of existence of
uniformly antisymmetric function
$f\colon\real\to\real$ with finite range remains open.
\end{abstract}

A function $f\colon\real\to\real$ is said to be
{\em uniformly antisymmetric\/} \cite{CL:Unif}
(or {\em nowhere weakly symmetrically continuous\/} \cite{thomson})
provided for every $x\in\real$ the limit
$\lim_{n\to\infty}(f(x+s_n)-f(x-s_n))$ equals $0$
for no sequence $\{s_n\}_{n<\omega}$ converging to $0$.
Uniformly antisymmetric functions have been studied
by Kostyrko~\cite{Kostyrko}, Ciesielski and  Larson~\cite{CL:Unif},
Komj\'ath and Shelah~\cite{KomShel},
and Ciesielski~\cite{CL:Unif1,CL:Unif2}.
(A connection of some of these results to the
paradoxical decompositions of the Euclidean
space $\real^n$ is described in Ciesielski \cite{45:SumAndDiff}.)
In particular in~\cite{CL:Unif} the authors constructed a
uniformly antisymmetric function $f\colon\real\to\N$
and noticed that the existence of a uniformly antisymmetric function
cannot be proved without an essential use of the axiom of choice.

The terminology and
notation used in this note is standard and follows~\cite{CiBook}.
In particular for a set $X$ we will write $|X|$ for its cardinality
and $\P(X)$ for its power set.
Also $2^\omega$ will stand for the set of all functions from
$\omega=\{0,1,2,\ldots\}$ into $2=\{0,1\}$.
We consider $2^\omega$ as ordered lexicographically.

\thm{thMain}{There exists a
function $f\colon\real\to\real$ with countable bounded range
such that for every $x\in\real$ there exists an $\ep_x>0$ with the property that
the set
\[
S_x=\{s\in\real\colon |f(x-s)-f(x+s)|<\ep_x\}
\]
is finite. In particular $f$ is uniformly antisymmetric.}

\proof
First notice that it is enough to find a compact zerodimensional
metric space $\la T,d\ra$ and a function $g$ from $\real$ into a countable
subset $T_0$ of $T$ such that
for every $x\in\real$ there is a $\delta_x>0$ for which the set
\[
\hat S_x=\{s\in\real\colon d(g(x-s),g(x+s))<\delta_x\}
\]
is finite.

To see this assume that such a function $g\colon\real\to T$ exists and take
a homeomorphic embedding $h$ of $T$ into $\real$.
We claim that $f=h\circ g\colon\real\to\real$
is as desired.
Indeed, $f[\real]=h[g[\real]]$ is countable, as it is a subset of a
countable set
$h[T_0]$, and it is bounded, since it is a subset of a compact set $h[T]$.
So take $x\in\real$ and $\delta_x>0$ for which $\hat S_x$ is finite.
Since $h^{-1}\colon h[T]\to T$ is uniformly continuous, we can find
an $\ep_x>0$ such that
\[
|y_1-y_2|<\ep_x\ \mbox{ implies }\ d(h^{-1}(y_1),h^{-1}(y_1))<\delta_x
\]
for every $y_1,y_2\in h[T]$. But for such a choice of $\ep_x$ we have
\[
S_x=\{s\in\real\colon |h(g(x-s))-h(g(x+s))|<\ep_x\}
\subset \hat S_x
\]
proving that $S_x$ is finite.

Thus, we proceed to construct a function $g$ described above.
The value of $g(x)$ will be defined with help of a representation of
$x$ in a Hamel basis, i.e., a linear basis of $\real$ over $\rational$.
For this we will use the following notation.
Let $\{y_\eta\colon \eta\in 2^\omega\}$ be a one-to-one enumeration of
a Hamel basis $\H$.
For every $x\in\real$ let $\sum_{\eta\in 2^\omega}q_{x,\eta}y_\eta$,
with $q_{x,\eta}\in\rational$ for $\eta\in 2^\omega$, be the unique
representation of $x$ in basis $\H$ and let
$w_x=\{\eta\in 2^\omega\colon q_{x,\eta}\neq 0\}$. Thus $w_x$ is finite and
\[
x=\sum_{\eta\in w_x}q_{x,\eta}y_\eta.
\]

The definition of the space $T$ is considerably more technical since it
reflects several different cases of  the
proof that the sets $\hat S_x$ are indeed finite.
To this end let
$\{q_j\colon j<\omega\}$ be a one-to-one enumeration of $\rational$ with
$q_0=0$. For $i<\omega$ let $\P_i=\P(\{q_j\colon j<i\})$ and
put $P_i=\P(2^i\times\{0,1\}\times\P_i\times\P_i)$.
Note that each $P_i$ is finite, so
$T=\prod_{i<\omega}P_i$, considered as the standard product
of discrete spaces, is compact zerodimensional.
We equip $T$ with a distance function $d$ defined
between different $s,t\in T$ by
\[
d(s,t)=2^{-\min\{i<\omega\colon s(i)\neq t(i)\}}
\]
and let
\[
T_0=\left\{t\in T\colon (\exists n<\omega)(\forall i\geq n)\
t(i)=\emptyset\right\}.
\]
Clearly $T_0$ is countable.

Now we are ready to define $g\colon\real\to T_0\subset T$.
For this, however, we will need few more definitions.
For $x\in\real$, $q\in\rational$, $i<\omega$, and $\zeta\in 2^i$ such that
$\zeta\in\{(\eta\restriction i)\colon \eta\in w_x\}$ we define:
\begin{itemize}
\item $p(i)\in\{0,1\}$ as the parity of $i$, i.e., $p(i)=i$ mod $2$;
\item $k_i(q)=\{q_j\in\rational\colon q_j<q\ \&\ j<i\}\in\P_i$;
\item $\eta(x,\zeta)$ to be the minimum of
$\{\eta\in w_x\colon\zeta\subset\eta\}$ (in the lexicographical order);
\item $\xi(x,\zeta)$ to be the minimum
of $\{\eta\in w_x\colon\zeta\subset\eta\}\setminus\{\eta(x,\zeta)\}$
provided
$|\{\eta\in w_x\colon\zeta\subset\eta\}|\neq 1$;
otherwise we put $\xi(x,\zeta)=\eta(x,\zeta)$;
\item $n_x<\omega$ to be the smallest number $n>0$ such that
\begin{description}
\item{(i)} $\eta\restriction n\neq \xi\restriction n$ for any different
$\eta,\xi\in w_x$, and

\item{(ii)} $q_{x,\eta}\in\{q_j\colon j<n\}$ for every
$\eta\in w_x$.
\end{description}
\end{itemize}
Consider the function $g\colon\real\to T_0$ defined as follows.
For every $x\in\real$ and $i<\omega$ we define $g(x)(i)\in P_i$ as
\[
\left\{\la\zeta,p(|\{\eta\in w_x\colon\zeta\subset\eta\}|),
k_i(q_{x,\eta(x,\zeta)}),k_i(q_{x,\xi(x,\zeta)})\ra\colon
\zeta\in\{(\eta\restriction i)\colon \eta\in w_x\}\right\}
\]
provided $i\leq n_x$
and we put $g(x)(i)=\emptyset$ for $n_x<i<\omega$.
In the argument below the key role will be played by
the function $k_i$ in general, and the coordinate
$k_i(q_{x,\eta(x,\zeta)})$ in particular.

The key step in the proof that $g$ has the desired property
is that for every $x\in\real$ and $s\neq 0$
\begin{equation}\label{eq1}
\mbox{if \ \ $n_x\leq\max\{n_{x-s},n_{x+s}\}$ \ \ then }\ \
g(x-s)(n_x)\neq g(x+s)(n_x).
\end{equation}

To see (\ref{eq1}) assume that $n_x\leq n_{x+s}$.
If $n_{x-s}<n_x$ then $g(x-s)(n_x)=\emptyset\neq g(x+s)(n_x)$, where
$g(x+s)(n_x)\neq\e$ since  $w_{x+s}\neq\e$ as
$n_{x-s}<n_x\leq n_{x+s}$ implies $x+s\neq 0$.
Thus, we can assume that $n_x\leq\min\{n_{x-s},n_{x+s}\}$.
Take an $\hat\eta\in w_{x-s}\cup w_{x+s}$ such that
$q_{x-s,\hat\eta}\neq q_{x+s,\hat\eta}$
and let $\zeta=\hat\eta\restriction n_x$.
Note that, by the definition of $n_x$,
the set $S=\{\eta\in w_x\colon\zeta\subset\eta\}$
has at most one element.

If $S=\emptyset$ then
$\{\eta\in w_{x-s}\colon\zeta\subset\eta\}=
\{\eta\in w_{x+s}\colon\zeta\subset\eta\}\neq\emptyset$
and so $\eta(x-s,\zeta)=\eta(x+s,\zeta)\notin w_x$
while $q_{x-s,\eta(x-s,\zeta)}+q_{x+s,\eta(x+s,\zeta)}=0$.
Thus $q_0=0$ separates
$q_{x-s,\eta(x-s,\zeta)}$ and $q_{x+s,\eta(x+s,\zeta)}$
implying that
$k_{n_x}(q_{x-s,\eta(x-s,\zeta)})\neq k_{n_x}(q_{x+s,\eta(x+s,\zeta)})$.
Therefore $g(x-s)(n_x)\neq g(x+s)(n_x)$.

So, assume that $S\neq\emptyset$ and let $\eta'$
be the only element of $S$.
Then $\eta'\in w_{x-s}\cup  w_{x+s}$.
If $\eta'$ belongs to precisely one of the sets
$w_{x+s}$ and $w_{x-s}$, say $w_{x+s}$,
then
$\{\eta\in w_{x+s}\colon\zeta\subset\eta\}=
\{\eta\in w_{x-s}\colon\zeta\subset\eta\}\cup\{\eta'\}$.
In particular,
$p(|\{\eta\in w_{x+s}\colon\zeta\subset\eta\}|)\neq
p(|\{\eta\in w_{x-s}\colon\zeta\subset\eta\}|)$
implying that $g(x-s)(n_x)\neq g(x+s)(n_x)$.

So, we can assume that $\eta'\in w_{x-s}\cap  w_{x+s}$.
Then $\{\eta\in w_{x-s}\colon\zeta\subset\eta\}=
\{\eta\in w_{x+s}\colon\zeta\subset\eta\}$
and $\eta(x-s,\zeta)=\eta(x+s,\zeta)$.
We will consider three cases.

\medskip

{\sc Case 1:} $\eta'\neq\eta(x-s,\zeta)=\eta(x+s,\zeta)$.
Then
$q_{x-s,\eta(x-s,\zeta)}+q_{x+s,\eta(x+s,\zeta)}=0$,
so $q_0=0$ separates
$q_{x-s,\eta(x-s,\zeta)}$ and $q_{x+s,\eta(x+s,\zeta)}$.
Thus
$k_{n_x}(q_{x-s,\eta(x-s,\zeta)})\neq k_{n_x}(q_{x+s,\eta(x+s,\zeta)})$
and $g(x-s)(n_x)\neq g(x+s)(n_x)$.

\medskip

{\sc Case 2:} $\eta'=\eta(x-s,\zeta)=\eta(x+s,\zeta)$ and
$q_{x-s,\eta(x-s,\zeta)}\neq q_{x+s,\eta(x+s,\zeta)}$.
Then
$q_{x-s,\eta(x-s,\zeta)}+q_{x+s,\eta(x+s,\zeta)}=2 q_{x,\eta'}$
and, by the definition of $n_x$, $q_{x,\eta'}\in\{q_j\colon j<n_x\}$.
Since $q_{x,\eta'}$ separates
$q_{x-s,\eta(x-s,\zeta)}$ and $q_{x+s,\eta(x+s,\zeta)}$
we conclude that
$k_{n_x}(q_{x-s,\eta(x-s,\zeta)})\neq k_{n_x}(q_{x+s,\eta(x+s,\zeta)})$
and $g(x-s)(n_x)\neq g(x+s)(n_x)$.

\medskip

{\sc Case 3:} $\eta'=\eta(x-s,\zeta)=\eta(x+s,\zeta)$ and
$q_{x-s,\eta(x-s,\zeta)}=q_{x+s,\eta(x+s,\zeta)}$.
Then
$Z=\{\eta\in w_{x-s}\colon\zeta\subset\eta\}\setminus\{\eta(x-s,\zeta)\}=
\{\eta\in w_{x+s}\colon\zeta\subset\eta\}\setminus\{\eta(x+s,\zeta)\}$
is non-empty, since it contains $\hat\eta$, and so
$\xi(x-s,\zeta)=\xi(x+s,\zeta)\notin w_x$.
Therefore, as in Case 1,
$q_{x-s,\xi(x-s,\zeta)}+q_{x+s,\xi(x+s,\zeta)}=0$,
so $q_0=0$ separates
$q_{x-s,\xi(x-s,\zeta)}$ and $q_{x+s,\xi(x+s,\zeta)}$.
Thus
$k_{n_x}(q_{x-s,\xi(x-s,\zeta)})\neq k_{n_x}(q_{x+s,\xi(x+s,\zeta)})$
and $g(x-s)(n_x)\neq g(x+s)(n_x)$.

This finishes the proof of (\ref{eq1}).

\medskip

Next, for every $x\in\real$ put $\delta_x=2^{-n_x}$.
To finish the proof of the theorem it is enough to show that
every $\hat S_x$ defined for such a choice of $\delta_x$
is a subset of a finite set
\[
Z_x=\{s\in\real\colon w_{x+s}\subset w_x\ \&\ n_{x+s}< n_x\}=
\left\{\sum_{\eta\in w_x}p_{\eta}y_\eta\colon
p_\eta\in\{q_j\colon j< n_x\}\right\}.
\]

Indeed, take an $s\in \hat S_x$. Then,
by (\ref{eq1}) and the definition of the distance function~$d$, we have
$\max\{n_{x-s},n_{x+s}\}<n_x$.
Notice also that 
if $n_{x-s}\neq n_{x+s}$, say $n_{x-s}<n_{x+s}$,
then $g(x-s)(n_{x+s})=\e\neq g(x+s)(n_{x+s})$
implying that 
$d(g(x+s),g(x-s))\geq 2^{-n_{x+s}}>2^{-n_x}=\delta_x$,
which contradicts $s\in \hat S_x$.
So, we have $n_{x-s}=n_{x+s}$.
To prove that $s\in Z_x$ it is enough to show that $w_{x+s}\subset w_x$.
But if it is not the case then there exists an $\eta\in w_{x+s}\setminus w_x$.
Moreover,
$q_{x+s,\eta}=-q_{x-s,\eta}\neq 0$
and
$\eta=\eta(x+s,\zeta)=\eta(x-s,\zeta)$,
where $\zeta=\eta\restriction n_{x+s}$.
In particular,
$q_0=0$ separates
$q_{x+s,\eta(x+s,\zeta)}$ and $q_{x-s,\eta(x-s,\zeta)}$.
Therefore
$k_{n_{x+s}}(q_{x-s,\eta(x-s,\zeta)})\neq k_{n_{x+s}}(q_{x+s,\eta(x+s,\zeta)})$
and $g(x-s)(n_{x+s})\neq g(x+s)(n_{x+s})$.
So
$d(g(x+s),g(x-s))\geq 2^{-n_{x+s}}>2^{-n_x}=\delta_x$
again contradicting $s\in \hat S_x$.
Thus, $w_{x+s}\subset w_x$ and $s\in Z_x$.

This finishes the proof of the theorem.

\end{document}